\documentclass[11pt,a4paper]{article}
\usepackage{t1enc}
\usepackage[latin1]{inputenc}
\usepackage[english]{babel}
\usepackage{times}
\usepackage{amsmath}
\usepackage{amssymb}
\usepackage{amsthm}
\usepackage{a4wide}
\usepackage{graphicx}        

\newtheorem{remark}{Remark}[section]

\newtheorem{definition}{Definition}[section]

\newcommand{\D}{\ {\rm d}}

\newcommand{\dualprod}[2]{ \langle #1 , #2 \rangle}
\newcommand{\R}{\mathbf{R}}

\title{Integral approach to sensitive singular perturbations}

\author{N. Meunier\and E. Sanchez-Palencia}
\author{{Nicolas Meunier}\footnote{Universit\'{e} Paris Descartes, MAP5, 45-47 Rue des Saints P\`{e}res
75006 Paris,
France,
Nicolas.Meunier@math-info.univ-paris5.fr}
{ and E. Sanchez-Palencia}\footnote{University of Pierre et Marie Curie and
CNRS, Institut Jean Le Rond D'Alembert, 4, place Jussieu, 75252
Paris, France, sanchez@lmm.jussieu.fr}}

\date{}

\begin{document}

\maketitle

\newcommand{\oH}{{\mathaccent'27 H}}


\begin{abstract}
We consider singular perturbation elliptic problems depending on a
parameter $\varepsilon$ such that, for $\varepsilon=0$ the
boundary conditions are not adapted to the equation (they do not
satisfy the Shapiro - Lopatinskii condition). The limit only holds
in very abstract spaces out of distribution theory involving
complexification and non-local phenomena. We give a very
elementary model problem showing the main features of the limit
process, as well as a heuristic integral procedure for obtaining a
description of the solutions for small $\varepsilon$. Such kind of
problems appear in thin shell theory when the middle surface is
elliptic and the shell is fixed by a part of the boundary and free
by the rest.
\end{abstract}
\section{\label{S1}Introduction}
\noindent

The main purpose of this paper is to give general ideas on a kind
of singular perturbations arising in thin shell theory when the
middle surface is elliptic and the shell is fixed by a part of the
boundary and free by the rest as well as an integral heuristic procedure reducing them to simpler problems. The system depends drastically on
the parameter $\varepsilon$ equal to the relative thickness of the
shell. It appears that the "limit problem" for $\varepsilon=0$ is
highly ill-posed. Indeed, the boundary conditions on the free
boundary are not "adapted" to the system of equations; they do not
satisfy the Shapiro - Lopatinskii (SL hereafter) condition.
Roughly speaking, this amounts to some kind of "transparency" of
the boundary conditions, which allow some kind of locally
indeterminate oscillations along the boundary, exponentially
decreasing inside the domain. This pathological behavior is only
concerned with $\varepsilon=0$. In fact, for $\varepsilon>0$ the
problem is "classical". When $\varepsilon$ is positive but small,
the "determinacy" of the oscillations only holds with the help of
boundary conditions on other boundaries, as well as the small
terms coming from $\varepsilon>0$.

In such kind of situations, the limit problem has no solution
within classical theory of partial differential equations, which
is uses distribution theory. It is sometimes possible to prove the
convergence of the solutions $u^{\varepsilon}$ towards some limit
$u^{0}$, but this "limit solution" and the topology of the
convergence are concerned with abstract spaces not included in the
distribution space.

After recalling the SL condition (section \ref{S2}), we give in
section \ref{S3} a very simple example of such a perturbation problem.
The geometry of the domain (an infinite strip) allows explicit
treatment by Fourier transform in the longitudinal direction. The
inverse Fourier transform within distribution theory is only
possible for $\varepsilon>0$, whereas for $\varepsilon=0$ it is
only possible in the framework of analytic functionals (higly
singular and not enjoying localization properties). This example
shows the prominent role of components with high frequency; for
small $\varepsilon$, the "smooth parts" (i. e. with small $|\xi|$)
of the solutions may be neglected with respect to  "singular ones"
(i. e. with large $|\xi|$). We also recall an example of elliptic
Cauchy problem (in fact Hadamard's counter-example) which exhibits
some relation with the limit problem.

In section \ref{S4}, we report the heuristic procedure of \cite{EgMeSa}. In this latter article,  we addressed a more complicated problem including
a variational structure, somewhat analogous to the shell problem,
but simpler, as concerning an equation instead of a system. It is
shown that the limit problem contains in particular an elliptic
Cauchy problem. This problem was handled in both a rigorous (very
abstract) framework and using a heuristic procedure for exhibiting
the structure of the solutions with very small $\varepsilon$. The
reasons why the solution goes out of the distibution space as
$\varepsilon$ goes to $0$ are then evident. In section \ref{S4} we present a simplified version
of the heuristic procedure involving only the essential facts of
the approximation, which are very much analogous to the method of
construction of a parametrix in elliptic problems \cite{Tay}, \cite{EgSc97}:

-Only principal (with higher differentiation order) terms are
taken into account.

-Locally, the coefficients are considered to be constant, their
values being frozen at the corresponding points.

-After Fourier transform ($x \to \xi$), terms with small
$\xi$ are neglected with respect to those with larger $\xi$ (which
amounts to taking into account singular parts of the solutions
while neglecting smoother ones). We note that this approximation,
aside with the two previous ones, lead to some kind of "local
Fourier transform" which we shall use freely in the sequel.

Another important ingredient of the heuristics is a previous
drastic restriction of the space where the variational problem is
handled. In order to search for the minimum of energy, we only
take into account functions such that the energy of the limit
problem is very small. This is done using a boundary layer method
within the previous approximations, i. e. for large $|\xi|$. This
leads to an approximate simpler formulation of the problem for
small $\varepsilon$, where  it is apparent that the
limit problem involves a smoothing operator and cannot have a
solution within distribution theory.

It should prove useful to give an example of a sequence of
functions converging to an analytical functional (but going out of
the distribution space, then leading to a "complexification"
phenomenon). It is known (\cite{Sc50}, \cite{GuCh64}) that (direct and inverse) Fourier
transform within distribution theory is only possible for
temperate distributions, not allowing functions with exponential
growth at infinity. The space of (direct or inverse) Fourier
transform of general distributions is noted $Z'$. It is a space of
analytical functionals: the corresponding test functions are
analytical rapidly decreasing functions, forming the space Z.

Let us consider the (non temperate) distribution (or function)
$\hat{u}(\xi)=\cosh ( \xi)$. The sequence 
$$\hat{u}^\lambda(\xi)=\left\{ \begin{split}& \cosh ( \xi) \textrm{ if }  |\xi  |< \lambda , \\ 
& 0 \textrm{ otherwise }
\end{split} 
\right.
$$
converges to $\hat{u} $  in the distribution sense as $\lambda $ goes to infinity. The inverse Fourier transforms $u^\lambda
(x)$ converge in $Z'$ to the analytical functionnal $u(x)$. The functions $\hat{u}^\lambda(\xi)$ are tempered and their inverse Fourier transforms are easily computed by hand. It appears that for large $\lambda $ 
$$u^\lambda (x) \approx \frac{e^\lambda }{2 \pi} \frac{1}{1+x^2} ( \cos ( \lambda x)+ x \sin ( \lambda x)).$$ 
It is then apparent that $u^\lambda (x)$ consists of a "nearly periodic" function with period tending to zero along with $1/\lambda $, multiplied by an "envelop" defined by $\frac{1}{1+x^2}$ and by the factor $\frac{e^\lambda }{2 \pi}$. Moreover, it should be noted that the amplitude is exponentially large with respect to the inverse of the period.  It is then
apparent that the limit is an "extremely singular" function as the "graph"
fills the entire plane. Moreover, it is clear (and may be
rigorously proved \cite{EgMeSa} that the sequence $u^{\lambda}$
goes out of the distribution space everywhere, not only in the
vicinity of $x=0$ as is suggested by the formal inverse Fourier
transform of $\cosh (\xi )= \Sigma _{n=0}^{+\infty} \frac{\xi ^{2n}}{(2n) !}$, which is 
$$u(x)= \Sigma _{n=0}^{+\infty} \frac{-i}{(2n)!}\delta ^{2n}(x),$$
apparently a singularity "of order
infinity" at the origin. This fact constitutes an example of the
property that elements of $Z'$ can only be tested with analytical
test functions, then not enjoying localization properties.

\begin{figure}
\includegraphics[width=.58\linewidth]{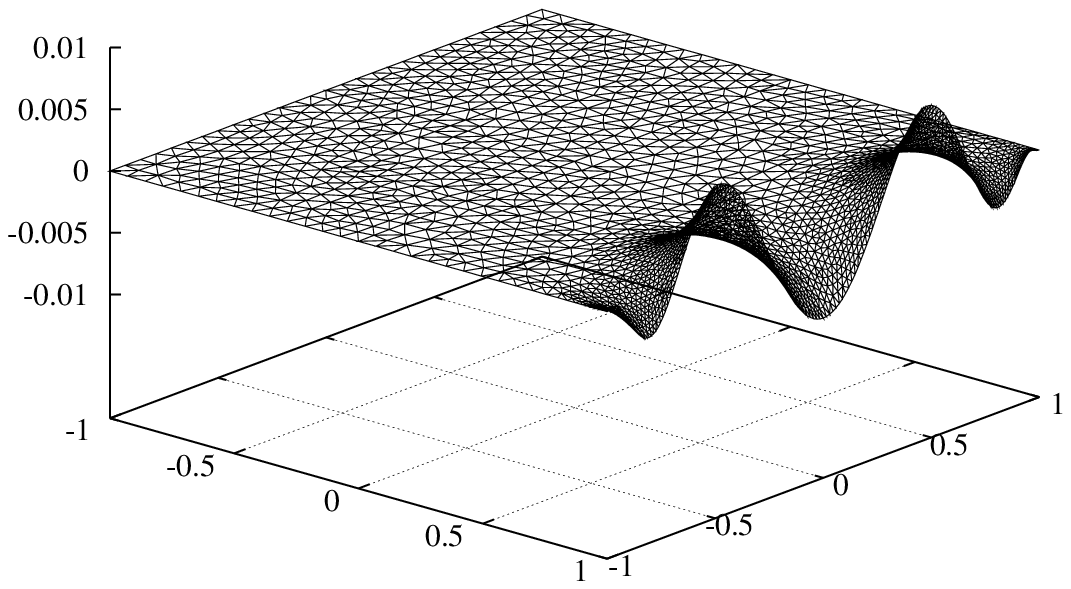}\,
\includegraphics[width=.58\linewidth]{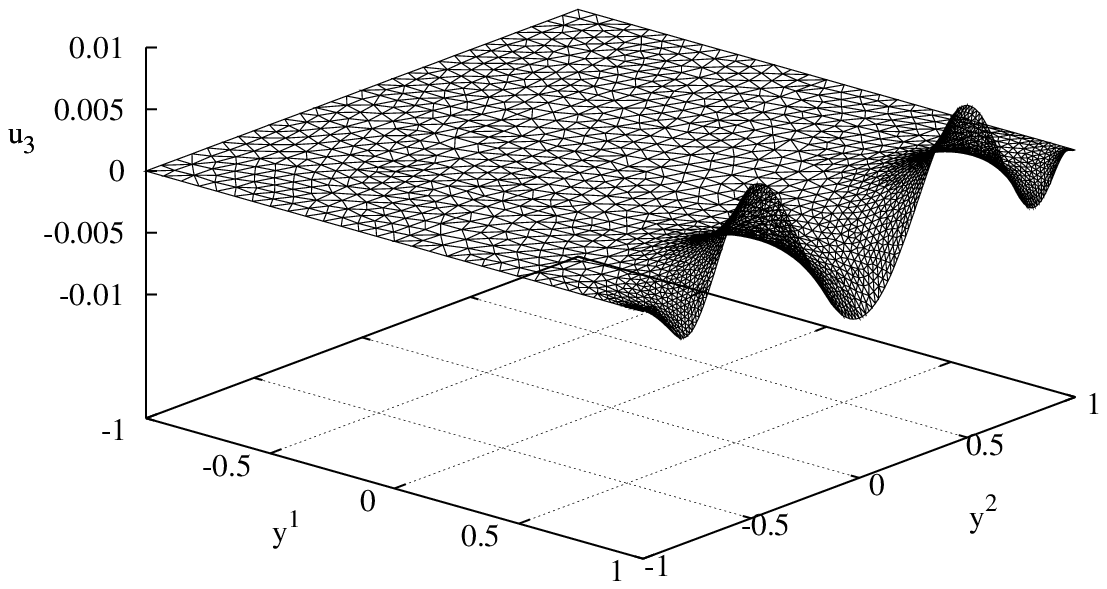}\,
\caption{Normal displacement for $\varepsilon = 10^{-3}$ on the left  and for $\varepsilon = 10^{-5}$ on the right }\label{F1}
\end{figure}

The motivation for studying that kind of problems comes from shell theory, see \cite{SaHuSa}, \cite{BeMiSa}. It appears that when the middle surface is elliptic (both principal curvatures have same sign) and is fixed by a part $\Gamma _0$ of the boundary and free by the rest $\Gamma _1$, the "limit problem" as the thickness $\varepsilon $ tends to zero is elliptic, with boundary conditions satisfying SL on $\Gamma _0$, and boundary not satisfying SL on $\Gamma_1$. Without going into details, which may be found in  \cite{MeunierSanchez},  \cite{Nantes}, \cite{EgMeSa} and \cite{EgMeSa09}, we show numerical computations taken from \cite{BeMiSa} of the normal displacement for $\varepsilon = 10^{-3}$ and $\varepsilon = 10^{-5}$ (figures \ref{F1} on the left and on the right respectively) when the shell is acted upon by a normal density of forces on a rectangular region of the plane of parameters. The most important feature is constituted by large oscillations nearby the free boundary $\Gamma _1$. It is apparent that, when passing from $\varepsilon = 10^{-3}$ to $\varepsilon = 10^{-5}$, the amplitude of the oscillations grows from 0.001 to 0.01. The singularities produced by the jump of the applied forces inside the domain is still apparent for $\varepsilon = 10^{-3}$, not for $\varepsilon =10^{-5}$, where only oscillations along the boundary are visible. Moreover, the number of such oscillations pass from nearly 3 for $\varepsilon =10^{-3}$ to nearly 5 for $\varepsilon =10^{-5}$ and is then nearly proportional to $\log (1/\varepsilon)$. We shall see that all these features agree with our theory.    

 \section{\label{S2} The Shapiro - Lopatinskii condition for boundary conditions of elliptic equations}
 \noindent
 
 In this section, we recall some properties of elliptic PDE, see \cite{AgDoNi59} and \cite{EgSc97} for more details.
   
 We consider a PDE of the form
 \begin{equation}
 P(x,\partial_{\alpha})u=f(x)
 \label{2.1}
\end{equation}
Where $x=(x_{1},x_{2})$ and $\partial_{\alpha}=
 \partial / \partial x _{\alpha}$, $\alpha=1,2$, and P is a polynomial of
degree $2m$ in $\partial_{\alpha}$. Let $P_{0}$ be the "principal
part", i. e., the terms of higher order. The equation is said to
be elliptic at $x$ if the homogeneous polynomial of degree $2m$ in
$\xi_{\alpha}$:
\begin{equation}
 P_{0}(x,-i \xi_{\alpha})=0
 \label{2.2}
\end{equation}
has no solution $\xi=(\xi_{1},\xi_{2})\neq (0,0)$ with real
$\xi_{\alpha}$. When the coefficients are real (this is the only
case that we shall consider) this implies that the degree is even
(this is the reason why we denoted it by $2m$). The left hand side
of (\ref{2.2}) is said to be the "principal symbol"; the "symbol" is
obtained in an analogous way taking the whole $P$ instead of the
principal part $P_{0}$. We note that replacing $\partial /
\partial x_{\alpha}$ by $-i\xi_{\alpha}$ in $P_{0}$ amounts to
taking formally the Fourier transform $x\rightarrow\xi$ for the
homogeneous equation with constant coefficients obtained by
discarding the lower order terms and freezing the coefficients at
$x$. Obviously, ellipticity on a domain $\Omega$ is defined as
elliptic at any $x \in \Omega$.

It is worthwhile mentioning that ellipticity amounts to non -
existence of "travelling waves" of the form
\begin{equation}
e^{-i\xi x}
 \label{m2.3}
\end{equation}
for the equation obtained after discarding lower order terms and
freezing coefficients. Here "travelling" amounts to "with real
$\xi$"; note that solutions as (\ref{m2.3}) with non real $\xi$ are
necessarily exponentially growing or decaying (in modulus) in some
direction. Moreover, when a solution of the form (\ref{m2.3}) exists (with
$\xi$ either real or not), it also exists for $c\xi$ with any $c$.
 In a heuristic framework,we may suppose that $|\xi|$ is very
 large; this justify to discard lower order terms (= of lower
 degree in $|\xi|$). In the same (heuristic) order of ideas,
 freezing the coefficients allows to consider "local solutions".
 This amounts to multiply the solutions by a "cutoff"
 function $\theta(x)$ or equivalently taking the convolution of the
 Fourier transform with $\hat{\vartheta}(\xi)$, which do not
 modify the behavior for large $\xi$. Microlocal analysis gives a
 rigorous sense to that heuristics. It then appears that local
 singularities of a solution $u$ (associated with behavior of the
 Fourier transform for large $|\xi|)$
  cannot occur in elliptic
 equations unless they are controlled by the (Fourier transform
 of the) right hand side $f$. This gives a "heuristic proof"
 of the classical property that local solutions of elliptic
 equations are rigorously associated with singularities of $f$.

 What happens with solutions near the boundary? Local Fourier
 transform is no longer possible, but, after rectification of the
 boundary in the neighborhood of a point, we may perform a
 tangential Fourier transform. If, for instance, the considered
 part of the boundary is on the axis $x_{1}$ and the domain is on the side $x_{2}>0$, taking only higher
 order terms and frozen coefficients, we have solutions of the
 form (\ref{m2.3}) with real $\xi_{1}$ (as coming from the Fourier
 transform) and non - real $\xi_{2}$. The dependence in $x_{2}$ is immediately obtained
  by solving an ODE with constant coefficients. Obviously, the solutions are
 exponentially growing or decreasing for $x_{2}>0$. As the
 coefficients are real, there are precisely $m$ (linearly independent) growing and $m$
 decreasing (in the case of multiple roots, dependence in $x_{2}$
 of the form $x_{2}e^{\lambda)x_{2}}$ and analogous also occur).
 Roughly speaking, there are solutions of the form:
 \begin{equation}
\sum_{k} C_{k } e^{-i\xi_{1} x_{1}}e^{\lambda_{k}x_{2}}
 \label{2.4}
\end{equation}
with real $\xi_{1}$ and $Re(\lambda) \neq 0$ (here $k$ is running
from $1$ to $2m$) . Boundary conditions on $x_{2}=0$ should
control solutions with $Re(\lambda) < 0$, i. e., exponentially
decreasing inside the domain, whereas exponentially growing ones
should be controlled "by the equation in the rest of the domain
and the boundary conditions on the other parts of the boundary".
In other words, "good boundary conditions" should determine,
(within our approximation of the half plane and frozen
coefficients) the solutions of the equation of the form (\ref{2.4}) with
$Re(\lambda) < 0$. Obviously, the number of such boundary
conditions is $m$. A set of $m$ boundary conditions enjoying the
above property is said to satisfy the Shapiro - Lopatinskii
condition. There are several equivalent specific definitions of
it. We shall mainly use the following one:

\begin{definition}
 Let $P$ be elliptic at a point $O$ of the boundary. A
set of $m$ boundary conditions
$B_{j}(x,\partial_{\alpha})=g_{j}(x)$, $j=1,...m$ is said to
satisfy the SL condition at $O$ when, after a local change to new
coordinates with origin at $O$ and axis $x_{1}$ tangent to the
boundary, taking only the higher order terms  and coefficients
frozen at $O$ in the equation and the boundary conditions, the
solutions of the form (\ref{2.4}) with $Re(\lambda) < 0$ obtained by
formal tangential Fourier transform are well defined by the
boundary conditions.
\end{definition}

\begin{remark}
The above definition should be understood in the sense of
formal solution for any given (real and non-zero) $\xi_{1}$. The
SL condition is not concerned with solutions in certain spaces. It
is purely algebraic, and concerns $m$ conditions imposed to the
$m$ (decreasing with $x_{2}$) linearly independent solutions of
the ODE obtained from $P_{0}$ by formal tangential Fourier
transform. This also amounts to saying that imposing the boundary
conditions equal to zero, the considered solutions must vanish. In
fact, the SL condition amounts to non-vanishing of a certain
determinant, and as so it is generically satisfied: conditions not
satisfying it are rarely encountered. In particular, in "well-behaved problems", when coerciveness on appropriate spaces is proved, the SL condition is not usually checked. It should also be noted that
the SL condition is independent of a change of variables,
and, in most cases, the change is trivial. On the other hand,
there are also definitions of the SL condition without change of
variables. Last, it should also be noted that the SL condition has
nothing to do with lower order terms and the right hand side of
the boundary conditions (as ellipticity is only concerned with the
principal symbol); it is merely a condition of adequation of the
principal part of the boundary operators to the principal part of
the equation.
\end{remark}

Let us consider, as an exercise, examples for the laplacian.
\begin{equation}
P = -\partial_{1}^{2}-\partial_{2}^{2}
 \label{2.5}
\end{equation}

The principal symbol is $\xi_{1}^{2}+\xi_{2}^{2}$ so that the
equation is elliptic of order $2$, then $m=1$. "Good boundary
conditions" are in number of 1.

Let us try the boundary condition (Dirichlet):
\begin{equation}
u = 0.
 \label{2.6}
\end{equation}
Taking any point of the boundary and $(x_{1},x_{2})$ with origin
at that point, tangent and normal to the boundary respectively,
the equation is the same as in the initial variables, and formal
tangential Fourier transform gives
\begin{equation}
 (\xi_{1}^{2}- \partial_{2}^{2}) \hat{u} (\xi _{1},x_{2}) = 0
 \label{2.7}
 \end{equation}

and the solutions are
\begin{equation}
 \hat{u}(\xi_{1},x_{2}) = C_{1}(\xi_{1}) e^{ |\xi_{1}|x_{2}}+
 C_{2 }(\xi_{1}) e^{|\xi_{1}|x_{2}}.
 \label{2.8}
\end{equation}
Taking only the exponentially decreasing for $x_{2}>0$ we only
have
\begin{equation}
  \hat{u}(\xi_{1},x_{2}) = C_{1 }(\xi_{1})e^{- |\xi_{1}|x_{2}}.
 \label{2.9}
\end{equation}
Now, imposing the "tangential Fourier transform" of
(\ref{2.6}):
\begin{equation}
  \hat{u}(\xi_{1},0) = 0,
 \label{2.10}
\end{equation}
we see that it vanish identically. Then, the Dirichlet boundary
condition satisfies the SL condition for the laplacian.

The case of the Neumann  boundary condition for the laplacian

\begin{equation}
\frac{\partial u}{\partial n }= 0.
 \label{2.11}
\end{equation}
is analogous. (Note also that the Fourier condition $(\frac{\partial u }{
\partial n })+ a u = g $ is
 the same, as only the higher order terms are taken in
 consideration). Proceeding as before, we have, instead of (\ref{2.10}):
\begin{equation}
  \partial _{2} \hat{u}(\xi _{1},0) = - |\xi_{1}| C_{1 }(\xi_{1})=0,
 \label{2.12}
\end{equation}
which also gives $C_{1 }(\xi_{1}) = 0$ and then $\hat{u} = 0$. Thus, (\ref{2.10}) satisfies SL for (\ref{2.5}).

Oppositely, the boundary condition:
\begin{equation}
 (\partial_{s} - i \partial_{n}) u = 0,
 \label{2.13}
\end{equation}
where $s$ and $n$ denote the arc of the boundary and the normal,
does not satisfy the SL condition for the laplacian. Indeed,
taking the new local axes, $s$ and $n$ become $x_{1}$ and $x_{2}$,
and after tangential Fourier transform:
\begin{equation}
 (-i \xi_{1} - i \partial_{2})\hat{u)}(\xi_{1},0) = 0,
 \label{2.14}
\end{equation}
which, applied to (\ref{2.9}) becomes:
\begin{equation}
 (-i \xi_{1} +i |\xi_{1}|)C_{1 }(\xi_{1})=0.
 \label{2.15}
\end{equation}
we then see that $C_{1 }(\xi_{1})$ vanishes for negative
$\xi_{1}$, but is arbitrary for positive $\xi_{1}$. In fact , the boundary condition
(\ref{2.13}) is "transparent" for solutions of the form
(\ref{2.9}) with positive $\xi_{1}$.

\begin{remark}
As it is apparent in the last example, when the SL
condition is not satisfied, there is some kind of "local
non-uniqueness", where "local" recalls that only higher order
terms are taken in consideration, and the coefficients are frozen
at the considered point of the boundary.
\end{remark}

The SL condition appears as some previous condition for solving
elliptic problems. It is apparent that some pathology is involved
at points of the boundary where it is not satisfied.

Let us mention, before closing this section, that the boundary
conditions may be different on different parts of the boundary
specially on different connected components of it (when there are
points of junction of the various regions, usually singularities
appear at that points).

\section{\label{S3}An explicit perturbation problem where the SL condition is not satisfied on a part of the boundary of the limit problem}
\noindent

Let $\Omega$ be the strip $(-\infty,+\infty) \times (0,1)$ of the
($x,y)$ plan. We denote by $\Gamma_{0}$ and $\Gamma_{1}$ the
boundaries $y=0$ and $y=1$ respectively. We then consider the
boundary value problem depending on the parameter $\varepsilon$:
 \begin{equation}
 \left\{
 \begin{array}{l}
\triangle u^{\varepsilon} =0 \textrm{ on }\Omega\\
u^{\varepsilon} =0 \textrm{ on } \Gamma_{0}\\
\partial_{x}u+ (i + \varepsilon^{2})\partial_{y}u = \varphi \textrm{ on  }\Gamma_{1}
\end {array}
\right. \label{3.1}
\end{equation}
where $\varphi$ is the data of the problem. It is a given function
of $x$, that we shall suppose sufficiently smooth, tending to $0$
at infinity. We shall solve it by $x\rightarrow \xi$ Fourier
transform; it is easily seen that we also have automatically $u
\rightarrow 0$ for $x\rightarrow \infty$, which may be added to
the boundary conditions.

The boundary condition on $\Gamma_{0}$ is the Dirichlet one, which
satisfies SL for the laplacian. Oppositely, the boundary condition
on $\Gamma_{1}$ satisfies it for $\varepsilon >0$ (this is easily
checked), not at the limit $\varepsilon =0$ (see the end of the previous section). The problem is to
solve for $\varepsilon >0$ and to study the behavior for
$\varepsilon $ going to zero.

Denoting by $\hat{}$ the $x\rightarrow\xi$ Fourier transform, $
\hat{u} ^{\varepsilon}$ is defined on the same $\Omega$ domain, but of the
$(\xi,y)$ plane. The solutions of the (transform of) equation and
the boundary condition on $\Gamma_{0}$ are of the form
\begin{equation}
\hat{u} ^{\varepsilon}(\xi,y) = \alpha(\xi)\sinh(\xi y)
  \label{3.2}
\end{equation}
where $\alpha$ denotes an unknown function to be determined with
the boundary condition on $\Gamma_{1}$. It will prove useful to
write the solution under the form
\begin{equation}
\hat{u} ^{\varepsilon}(\xi,y) = \hat{\beta}^{\varepsilon}(\xi) \frac{\sinh(\xi
y)}{\sinh(\xi)}
  \label{3.3}
\end{equation}
for the new unknown $\hat{\beta}^{\varepsilon }(\xi)$, which is the transform of
the trace $u^{\varepsilon}(x,0)$. Imposing the Fourier transform
of the boundary condition on $\Gamma_{1}$ we have:

\begin{equation}
-i \xi \hat{\beta}^{\varepsilon}(\xi) + (i + \varepsilon^{2})\frac{\cosh
(\xi)}{\sinh (\xi)}\hat{\beta}^{\varepsilon} (\xi) \xi = \hat{\varphi}(\xi).
  \label{3.4}
\end{equation}
 So that:
\begin{equation}
\hat{\beta}^{\varepsilon}(\xi) =\frac{\hat{\varphi}(\xi)}{-i \xi \Big(1
-\coth (\xi) \Big)+ \varepsilon^{2}\xi \coth (\xi) } .
  \label{3.5}
\end{equation}
In order to study this function, we should keep in mind that the
expression $(1 - \coth (\xi))$ decays for $\xi\rightarrow+\infty$ as
$2e^{-2\xi}$. Then, at the limit $\varepsilon = 0$ we have
\begin{equation}
\hat{\beta}^{0}(\xi) =\frac{\hat{\varphi}(\xi)}{-i \xi(1 - \coth (
\xi))}.
  \label{3.6}
\end{equation}
For $\xi\rightarrow+\infty$ this function behaves as
\begin{equation}
\hat{\beta}^{0}(\xi) \approx 2 \frac{\hat{\varphi}(\xi)}{-i
\xi}e^{2\xi}.
  \label{3.7}
\end{equation}
This shows (unless in the case of very special data $\varphi$ with very fast
decaying Fourier transform) that $\hat{\beta}^{0}(\xi)$ is not a
tempered distribution, and the inverse Fourier transform is an analytical function in $\mathcal{Z}'$. Nevertheless, for $\varepsilon >
0 $, $\hat{\beta}^{\varepsilon}(\xi)$ is "well-behaved" for
$\xi\rightarrow +\infty$ as
\begin{equation}
\hat{\beta}^{\varepsilon}(\xi) \approx \frac{\hat{\varphi}(\xi)}{
\xi \varepsilon^{2}}.
  \label{3.8}
\end{equation}
This specific behavior depends on that of $\frac{\hat{\varphi}}{ \xi}$,
so that in most cases will be decreasing, but multiplied by the
factor $\varepsilon ^{-2}$. When $\varepsilon >0$ (small but not $0$)  is fixed,  $\hat{\beta}^{\varepsilon}(\xi)$ is approximatively given by  (\ref{3.6}) for "finite" $\xi $ and by  (\ref{3.8}) for $\xi $ going to $+\infty$. It is easily seen that the  sup in
modulus of $|\hat{\beta}^{\varepsilon}(\xi)|$ is located in the
region where both terms in the denominator of the right hand side
of (\ref{3.5}) are of the same order (so that no one of them may
ble neglected). This gives
\begin{equation}
\xi = \mathcal{O}(log(1/\varepsilon)).
  \label{3.9}
\end{equation}
 
It appears that $\hat{\beta}^{\varepsilon}(\xi)$ consists mainly of Fourier components which tend to infinity algebraically as $\varepsilon $ goes to zero with $\xi $ tending to infinite "slowly" as in  (\ref{3.9}). this is somewhat analogous to the example, given in the introduction, of a sequence of functions converging to an analytical functional. 

Coming back to (\ref{3.3}), the main properties of the behavior of $u^{\varepsilon}(x,1)$
 may be thrown:

- The trace $u^{\varepsilon}(x,1)= \beta^{\varepsilon}(x)$ on the
boundary $\Gamma_{1}$ which bears the "pathological boundary
condition" mainly consists of large oscillations with wave length
$1/\log(1/\varepsilon)$ (which tends to $0$ very slowly as
$\varepsilon \to 0$). The amplitude of that oscillations grows
nearly as $\varepsilon^{-2}$. The limit $\varepsilon \to 0$ does
not exist in distribution theory; it constitutes a
complexification process.

- Out of the trace on $\Gamma_{1}$, (i. e. for $0<y<1$), the
behavior is analogous, but of lower amplitude, which is
exponentially decreasing going away of $\Gamma_{1}$. We recover
properties of the non-uniqueness associated with the failed SL
condition.

Before concluding this section, we would like to show some analogy
between the previous limit problem and the Cauchy elliptic
problem, which is a classical example of ill-posed problem,
without solution in general.

We consider the same domain $\Omega$ as before, but we now impose
two boundary conditions on $\Gamma_{0}$ and no condition on
$\Gamma_{1}$. Namely
 \begin{equation}
 \left\{
 \begin{array}{l}
\triangle v =0 \textrm{ on } \Omega\\
v =\psi \textrm{ on } \Gamma_{0}\\
\partial_{y}v=0  \textrm{ on } \Gamma_{0}
\end {array}
\right. \label{3.10}
\end{equation}

Taking as above the $x\rightarrow\xi$ Fourier transform, it
follows immediately that
\begin{equation}
\hat{v}(\xi,y) = \hat{\psi}(\xi)\cosh(\xi y).
  \label{3.11}
\end{equation}
Where it is apparent that the behavior for $\xi\rightarrow\infty$
is exponentially growing (unless in the case when
$\hat{\psi}(\xi)$  decays faster than $e^{-|\xi|}$) so that it is
not tempered and the inverse Fourier transform does not exist
within distribution theory.

\section{\label{S4} A model variational sensitive singular perturbation,  \cite{EgMeSa}  }
\subsection{Setting of the problem}
\noindent

Let $\Omega $ be a two dimensional compact manifold with smooth
(of $C^{\infty }$ class) boundary $\partial \Omega = \Gamma _0
\cup \Gamma _1$ of the variable $x=(x_1,x_2)$, where $\Gamma _0 $
and $\Gamma _1$ are disjoint; they are one - dimensional compact
smooth manifolds without boundary, then diffeomorphic to the unit
circle. Let $a$ and $b$ be the bilinear forms given by:
\begin{eqnarray}\label{defa}
a(u,v)&=&\int _\Omega \triangle u \, \triangle v \D x, \\
b(u,v)&=&\int _\Omega \sum_{\alpha ,\beta =1}^{2} \partial_{\alpha
\beta} u \, \partial _{\alpha \beta} v \D x .\label{defb}
\end{eqnarray}
We consider the following variational problem (which has possibly
only a formal sense)
\begin{equation}\label{pbvar}
\left\{\begin{array}{l}\textrm{Find } u ^{\varepsilon }\in V
\textrm{ such that, }
\forall v \in V\\
a(u^{\varepsilon } ,v)+\varepsilon ^2 b(u^{\varepsilon }, v)
=\dualprod{f}{v},
\end{array}\right.
\end{equation}
where the space $V$ is the "energy space" with the essential
boundary conditions on $\Gamma _0$
\begin{equation}\label{defdeV}
V=\{v \in H^2(\Omega ); \ v_{| \Gamma _0}=\frac{\partial
v}{\partial n}_{| \Gamma _0} =0\},
\end{equation}
where $n,t$ denotes the normal and tangent unit vectors to the
boundary $\Gamma $ with the convention that the normal vector $n$
is inwards $\Omega $. It is easily checked that the bilinear form
$b$ is coercive on $V$. Moreover, we immediately obtain the
following result. For all $\varepsilon >0$ and for all $f $ in
$V'$, the variational problem (\ref{pbvar}) is of Lax-Milgram type
and it is a self-adjoint problem which has a coerciveness constant
larger than $c \varepsilon ^2$, with $c>0$.

The equation on $\Omega$ associated with problem (\ref{pbvar}) is:
\begin{equation}\label{pbfor}
(1+ \varepsilon ^2) \triangle ^2 u^{\varepsilon} =f \textrm{ on
}\Omega ,
\end{equation}
as both forms $a$ and $b$ give the laplacian. As for the boundary
conditions on $\Gamma_{0}$, they are "principal" i. e. they are
included in the definition of $V$, (5.4). As for conditions on
$\Gamma_{1}$, they are "natural", classically obtained from the
integrated terms by parts. Those coming from the form $b$ are
somewhat complicated; we shall not write them, as the problem with
$\varepsilon>0$ is classical. For $\varepsilon=0$ these conditions
(coming from form $a$) are: $ \triangle u= \frac{\partial \triangle u}{\partial n }=0,
\textrm{ on }\Gamma _1$.

As a matter of fact, the full limit boundary boundary value
problem is:
\begin{equation}\label{pbvar3bis}
\left\{\begin{array}{l}
 \triangle ^2 u^{0} =f \textrm{ on }\Omega \\
u= \frac{\partial u^{0}}{\partial n }=0, \textrm{ on }\Gamma _0\\
\triangle u ^{0}=0  \textrm{ on  } \Gamma _1 \\
- \frac{ \partial }{\partial n }  \triangle u ^{0} =0  \textrm{ on
} \Gamma _1 .
\end{array}\right.
\end{equation}

Let us check that the boundary conditions on $\Gamma_{1}$ (i. e;
the two last lines of (5.6)) do not satisfy the SL condition for
the elliptic operator $\triangle^{2}$.
   Indeed, proceeding as in sect. 2, by formal tangential Fourier
transform

\begin{equation}\label{bande1}
(-\xi _1 ^2+ \partial _2 ^2)^2\hat{u} =0.
\end{equation}
which yields that
\begin{equation}\label{bande1bis}
\hat{v}=  (A e^{-|\xi _1 | x_2 }+C x_2 e^{-|\xi _1 | x_2 })
\end{equation}
(as well as analogous terms with $+|\xi|$ instead of $-|\xi|$,
which are not taken into account as exponentially growing inwards
the domain). Here, according to SL theory, $x_{2}$ is the
coordinate normal to the boundary, after taking locally tangent
and normal axes, (which do not modify the equation
$\triangle^{2}$). The (tangential Fourier transform of the)
boundary conditions on $\Gamma_{1}$ are:
\begin{equation}
(-\xi _1 ^2+ \partial _2 ^2)\hat{u} =0
\end{equation}
and
\begin{equation}
\partial_{2}(-\xi _1 ^2+ \partial _2 ^2)\hat{u} =0.
\end{equation}
It is immediately seen that the previous solutions (\ref{bande1bis}) with
$C=0$ and any $A\neq0$ satisfy both conditions (note that its
laplacian vanishes everywhere, then it vanishes as well as its
normal derivative on the boundary). So, the SL condition is not
satisfied on $\Gamma_{1}$.

Before going on with our study, we note that the limit problem
(\ref{pbvar3bis}) implies an elliptic Cauchy problem for the auxiliary unknown
\begin{equation}
v^{0}= \triangle u^{0}.
\end{equation}
Indeed, system (\ref{pbvar3bis}) gives in particular:
\begin{equation}
\left\{\begin{array}{l}
 \triangle v^{0} =f \textrm{ on }\Omega \\
 v^{0}=0  \textrm{ on  } \Gamma _1 \\
- \frac{ \partial  v^{0}}{\partial n }   =0  \textrm{ on
} \Gamma _1 .
\end{array}\right.
\end{equation}
which is precisely the Cauchy problem for the laplacian.

As mentioned in section  \ref{S3}, this is a classical ill - posed problem,
and the solution does not exist in general. Oppositely, uniqueness
of the solution holds true (uniqueness theorem of Holmgren and analogous, see for instance \cite{CoHi62}).

\subsection{The heuristic integral approach}
\noindent

The aim of this section is the construction, in a heuristic way,
of an approximate description of the solutions $u^{\varepsilon}$
of the model problem in the previous section for small values of
$\varepsilon$.

From the general theory of singular perturbations of the form
(\ref{pbvar}), we know that our assumption
\begin{equation}\label{3.1}
a(v,v)^{1/2} \textrm{ defines a norm on }V,
\end{equation}
is crucial. Indeed, when it is not satisfied, the problem is said
to be "non inhibited". In such a case, it has a kernel which
contains non vanishing terms and then, it is easy to establish
that the asymptotic behaviour of the solution $u^{\varepsilon }$
of (\ref{pbvar}) is described by a variational problem in this
kernel. The previous fact is not surprising as soon as we consider
the following minimization problem, which is equivalent to
(\ref{pbvar}),
\begin{equation}\label{pbmin}
\left\{\begin{array}{l}\textrm{Minimize in }  V,\\
a(u^{\varepsilon } ,u^{\varepsilon })+\varepsilon ^2
b(u^{\varepsilon }, u^{\varepsilon }) -
2\dualprod{f}{u^{\varepsilon }}.
\end{array}\right.
\end{equation}
Indeed,  when
$\varepsilon $ goes to zero, the natural trend consists in
avoiding the $a$-energy which occurs with the factor $1$ and
leaving the $b$-energy which has a factor $\varepsilon ^2$.

Clearly, this is not possible when (\ref{3.1}) is satisfied since
the kernel reduces to the zero function. Nevertheless, in our case, $a(v,v)=0$ implies $\triangle v =0$ and, as $v \in V$, the traces of $v$ and $\frac{\partial v}{\partial n}$ vanish on $\Gamma _0$, so that   (\ref{3.1}) follows from
 the uniqueness theorem for the Cauchy problem. This uniqueness is classical, but the solution $u$ is unstable in the
sense that there can be "large $u$" in the $V$ norm (or in other
spaces) for "small $f$" in the $V ' $ norm (or in other spaces).
It then appears that the same reasoning shows that for small
values of $\varepsilon $, the solution $u^{\varepsilon }$ will
 be precisely among elements with small $a(u^{\varepsilon }, u^{\varepsilon })$, that is to say
with small $\triangle u^{\varepsilon }$ in $L^2$.

\subsection{The $\Gamma _0$ layer}
\noindent

Let us now build such functions $u^{\varepsilon } \in V$ with very small$\|\triangle u^{\varepsilon } \|_{L^2}$. The main
idea is to consider functions in a larger space than the space of
functions $v$ of $V$ such that $\triangle v=0$ (which only
contains the function $v=0$). The functions of this bigger space
will not satisfy the two boundary conditions on $\Gamma _0 $ that
are satisfied by any function of $V$.  Then we shall modify it in
a narrow boundary layer along $\Gamma _0$ in order to satisfy the
two boundary conditions with small value of $a$-energy.

More precisely, let us consider the vector space:
\begin{equation}\label{definitiondeg0}
G^0=\{v \in C^{\infty }(\overline{\Omega}), \ \triangle v=0
\textrm{ on } \Omega , \ v=0 \textrm{ on } \Gamma _0\}.
\end{equation}

\begin{remark}
We observe that every function of $G^0$ satisfies one of the
boundary conditions on $\Gamma _0$ which are satisfied by any
element of $V$. For simplicity, we have chosen $v=0$ on $\Gamma
_0$, but we could choose the other one $\frac{\partial v}{\partial
n}=0$ on $\Gamma _0$ as well. On the other hand, the regularity
assumption $C^{\infty }$ is slightly arbitrary. Since, we will
consider the completion of $G^0$ with respect to some norm, this
point is irrelevant.
\end{remark}

Obviously, as the Dirichlet problem for the laplacian on $\Omega $
is well posed in $C^{\infty }$, the space $G^0$ is isomorphic with
the space of traces on $\Gamma _1$:
\begin{equation}
\{w \in C^{\infty }(\Gamma _1)\}
\end{equation}
the isomorphism is obtained by solving the Dirichlet problem:
\begin{equation}\label{pbdirichlet}
\left\{\begin{array}{l}\triangle \tilde{w} =0\textrm{ on }  \Omega,\\
\tilde{w} =0\textrm{ on }  \Gamma _0,\\
\tilde{w} =w\textrm{ on }  \Gamma _1.
\end{array}\right.
\end{equation}

In the sequel, we shall consider indifferently the functions
$\tilde{w} $ on $\overline{\Omega }$ or their traces $w$ on
$\Gamma _1$.

In fact, the exact function $u^{\varepsilon}$ is a solution of
(\ref{pbfor}), which we are searching to describe approximatively
in order to define a space as small as possible (incorporating the
main features of the solution) to solve the minimization problem.
More precisely, according to our previous comments, we are
interested in the "most singular parts" of $u^{\varepsilon}$ in
the sense of the part corresponding to the high frequency Fourier
components. As we shall see in the sequel, it turns out that these
singular parts may be obtained by modification of the functions
$\tilde{w}$ on a boundary layer close to $\Gamma_{0}$; this layer
is narrower when the considered Fourier components are of higher
frequency; in fact, the layer only exists because we only consider
high frequencies. This allows to make an approximation which
consists in using locally curvilinear coordinates defined by the
arc of $\Gamma_{0}$ and the normal, and handling them as cartesian
coordinates. Clearly, this approximation is exact only on the very
$\Gamma_{0}$, but more and more precise as we approach of
$\Gamma_{0}$, i. e. as the considered frequencies grow.

Once the layer is constructed, we compute the $a$-energy of it, as
well as the $\varepsilon^{2}b$-energy of the (modified)
$\tilde{w}$ function, in order to consider the variational problem
(\ref{pbvar}) in the restricted space.

Let us first exhibit the local structure of the Fourier transform
of $\tilde{w}$ close to $\Gamma_{0}$. According to our general
considerations on the heuristic procedure, $\hat{w}$ may be considered (after multiplying by an appropriate cutoff function) of "small support" near a point $P_0$ of $\Gamma _0$. Taking local tangent and normal cartesian coordinates $y_1, y_2$, we have, within our approximation,
\begin{equation}\label{approx}
 \Big(\frac{\partial ^2}{\partial y_1^2}+\frac{\partial ^2}{\partial y_2^2}\Big) \tilde{w}  =0 \textrm{ on } \R \times (0,t),
\end{equation}
for some $t>0$.
Taking the tangential Fourier transform, we obtain:
\begin{equation}\label{solmodifie2}
\mathcal{F}(\tilde{w}_j) (\xi _1, y_2)=\lambda e^{|\xi _1 | y_2}
+\mu e^{-|\xi _1 | y_2}.
\end{equation}
It is worthwhile defining the local structure of $\hat{w}$ in the
vicinity of $\Gamma_{0}$ using the "Cauchy" data $\tilde{w}$ and
$\partial_{2}\tilde{w}$ on $\Gamma_{0}$ (note that the solution of
the Cauchy problem is unique, so that the Cauchy data determine
the solution). As $\hat{w}$ vanishes on $\Gamma_{0}$, the local
structure is then determined by $\partial_{2}\tilde{w}$ on
$\Gamma_{0}$. Taking the tangential Fourier transform this gives:
\begin{equation}\label{solmodifie5}
\mathcal{F}\Big(\tilde{w} _j\Big)(\xi _1,
y_2)=\mathcal{F}\Big(\frac{ \partial \tilde{w}_j }{\partial
y_2}_{| y_2=0}\Big) \frac{\sinh (|\xi _1 | y_2)}{|\xi _1 |}.
\end{equation}

We now proceed to the modification of $ \tilde{w}$ into
$\tilde{w}^ a$ in a narrow boundary layer of $\Gamma _0$ in order
to satisfy (always within our approximation) the equation coming
from (\ref{pbfor}) for small $\varepsilon $. Using considerations
similar to those leading to (\ref{approx}), this amounts to
\begin{equation}\label{eqmodifie1}
\Big( \frac{\partial ^2 }{\partial y_1 ^2}+\frac{\partial ^2
}{\partial y_2 ^2} \Big)^{(2)} \tilde{w}^a= 0 \textrm{ on } \R
\times (0,t).
\end{equation}
hence the tangential Fourier transform reads
\begin{equation}\label{eqmodifie2}
\Big( -|\xi _1 |^2+\frac{\partial ^2 }{\partial y_2 ^2}
\Big)^{(2)} \mathcal{F}(\tilde{w}^a) = 0.
\end{equation}

Consequently, $\mathcal{F}(\tilde{w}^a)$ should take the form
\begin{equation}\label{solmodifie}
\mathcal{F}(\tilde{w}^a) (\xi _1, y_2)=(\alpha + \gamma
y_2)e^{|\xi _1 | y_2} +(\beta + \delta y_2)e^{-|\xi _1 | y_2}.
\end{equation}

The four unknown constants should be determined by imposing that $
\tilde{w}^a$ and $ \partial_{2}\tilde{w}^a$ vanish for $y_{2}=0$
and the "matching condition" of the layer, i.e., out of the layer,
we want $\tilde{ w} ^a _j$ to match with the given function $
\tilde{ w}_j$. Since $|\xi _1 |>>1$, then $|\xi _1 | y_2 >>1$
means that $y_2 >>\frac{1}{|\xi _1|}$ (but we still impose that
$y_2$ is small in order to be in a narrow layer of $\Gamma _0$);
this is perfectly consistent, as we will only use the functions
for large $|\xi_{1}|$, hence the terms with coefficients $\beta $
and $\delta $ are "boundary layer terms" going to zero out of the
layer (i.e. for $|y_2 |>>\mathcal{O}\Big( \frac{1}{|\xi
_1|}\Big)$), see perhaps \cite{Ec79} or \cite{Il'in} for generalities on boundary layers and matching. This gives
\begin{equation}\label{solmodifie5}
\mathcal{F}\Big(\tilde{w} _j\Big)(\xi _1,
y_2)=\mathcal{F}\Big(\frac{ \partial \tilde{w}_j }{\partial
y_2}_{| y_2=0}\Big)( \frac{\sinh (|\xi _1 | y_2)}{|\xi _1 | }-
y_2e^{-|\xi _1 | y_2} ).
\end{equation}

This amounts to saying that the modification of the function
$\tilde{w}_{j}$ consists in adding to it the inverse Fourier
transform of
\begin{equation}\label{V1}
 \mathcal{F}\Big(\frac{
\partial \tilde{w} _j}{\partial y_2}_{| y_2=0}\Big)
 \Big(
 - y_2e^{-|\xi _1 | y_2}\Big).
\end{equation}

Defining on $\Gamma_{0}$ the family (with parameter $y_{2}$) of
pseudo-differential smoothing operators
$\delta\sigma(\varepsilon,D_{1},y_{2})$ with symbol:
\begin{equation}\label{V5}
 \delta \sigma (\varepsilon,\xi _1, y_2)= -y_2 e^{- |\xi _1 |
 y_2}h(\varepsilon,\xi , y_2),
\end{equation}
where $h$ is an irrelevant cutoff function avoiding low frequencies; it is equal to 1 for high frequencies (see \cite{EgMeSa} for details), we see that the modification of the function $ \tilde{w}$:
\begin{equation}\label{V6}
\delta \tilde{w} =  \tilde{w}^{a}- \tilde{w}
\end{equation}
is precisely the action of $\delta\sigma(\varepsilon,D_{1},y_{2})$
on $ \frac {\partial \tilde{w}_{j}}{\partial y_2}(y_{1},0)$:
\begin{equation}\label{V7}
\delta \tilde{w} =\delta\sigma(\varepsilon,D_{1},y_{2})\frac
{\partial \tilde{w}_{j}}{\partial y_2}(y_{1},0).
\end{equation}

Let us now compute the leading terms of the $a$-energy  of the
modified function $\tilde{w}^a$.

Let $\tilde{v}$ and $\tilde{w}$ be two elements in $G^{0}$ and
$\tilde{v}^{a}, \tilde{w}^{a}$ the corresponding elements modified
in the boundary layer. As the given $\tilde{v}$ and $\tilde{w}$
are harmonic in $\Omega $, the $a$-form is only concerned with the
modification terms $\delta\tilde{v}$ and $\delta\tilde{w}$. Then,
within our approximation, we have:
\begin{equation}\label{V8}
a( \tilde{v} ^a, \tilde{w}^a)=  \int _{\Gamma _0 } dy_1 \int _0
^{+\infty } \triangle (\delta \tilde{v})\overline{\triangle
(\delta \tilde{w}})\D y_2.
\end{equation}

To compute this expression, we first write $\tilde{v}$ and $\tilde{w}$ as sum of terms with "small support" (by multiplying by a partition of unity): $\tilde{v}=  \Sigma _{j}\tilde{v}_j$ and $\tilde{w}= \Sigma _{j} \tilde{w}_j$. Then, within our approximation, the integral is on the halfplane $\R \times (0, +\infty)$ of the variables $y_1,y_2$. Taking the tangential Fourier transform 
and using the Parceval-Plancherel theorem, we have
\begin{eqnarray*}\label{energiea2}
a(\tilde{v}^a, \tilde{w}^a)=  \Sigma _{j,k} \int _{-\infty } ^{+
\infty} \D \xi _1 \int _0 ^{+\infty } \Big(\frac{\D ^2}{\D y_2^2}-
\xi _1 ^2\Big)  \delta \sigma (\varepsilon,\xi, y_2)\mathcal{F}
\big(\frac{ \partial \tilde{v}_{j} }{\partial y_2}_{| y_2=0}\big)\times \nonumber \\
\overline{\Big(\frac{\D ^2}{\D y_2^2}- \xi _1 ^2\Big)  \delta
\sigma (\varepsilon,\xi, y_2) \mathcal{F} \big(\frac{ \partial
\tilde{w}_{k} }{\partial y_2}_{| y_2=0}\big)}\D y_2 .
\end{eqnarray*}
Hence, on account of (\ref{V5}) and integrating in $y_2$, this yields 
\begin{equation}\label{energiea1}
a(\tilde{v}^a, \tilde{w}^a)=\Sigma _{j, k}\int _{-\infty } ^{+
\infty} 2 |\xi _1| \frac{\partial \tilde{w}_{1,j} }{\partial
y_2}_{|y_2=0} \overline{\frac{\partial \tilde{w}_{2,k}}{\partial
y_2}_{|y_2=0}} h^{2}(\varepsilon,\xi , y_2) \D \xi _1.
\end{equation}

This expression (\ref{energiea1}) only depends on the traces
$\frac{\partial \tilde{v}_{j} }{\partial y_2}_{|y_2=0}(y_1)$ and
 $\frac{\partial \tilde{w}_{k} }{\partial y_2}_{|y_2=0}(y_1)$, which
are functions defined on $\Gamma _0$.

We now simplify this last expression using a sesquilinear form
involving pseudo-differential operators.

Indeed, denoting by $P(\frac{\partial }{\partial y_1}) $ the
pseudo-differential operator with symbol
\begin{equation}\label{pseudodiff}
P(\xi _1 )= (2 |\xi _1 |)^{1/2}h(\varepsilon,\xi , y_2),
\end{equation}
and summing over $j$ and $k$, we obtain that
\begin{equation}\label{energiea3bis}
a(\tilde{v}^a, \tilde{w}^a)=  \int _{\Gamma _0} P(\frac{\partial
}{\partial s}) \frac{\partial \tilde{v} }{\partial n}_{|\Gamma
_0}\overline{P(\frac{\partial }{\partial s})
 \frac{\partial \tilde{w} }{\partial n}_{|\Gamma _0}} \D s.
\end{equation}

\subsection{Taking account of the perturbation term $\varepsilon ^2b$.}\label{section2.4}
\noindent

We now consider the minimization problem (\ref{pbmin}) on $G^0$
instead of on $V$. Obviously, the $a$-energy should be computed
using formula (\ref{energiea3bis}). This modified problem
should involve the $a$-energy and the
$\varepsilon^{2}b$-energy. A natural space for handling it should
be the completion $G$ of $G^0$ with the norm:
\begin{equation}\label{normG}
\|v \| _{G} ^2 = \int _{\Gamma _0} \Big|P(\frac{\partial
}{\partial s}) \frac{\partial v }{\partial n}_{|\Gamma _0}\Big |^2
\D s + b(v,v).
\end{equation}

It is easily seen that $G$ is the space of the harmonic functions
of $H^{2}(\Omega)$ vanishing on $\Gamma_{0}$; according to
(\ref{pbdirichlet}) it may be identified with the space of traces
$H^{3/2}(\Gamma_{1})$.

It will prove useful to write another (asymptotically equivalent
for large $|\xi _1|$) definition of this problem. Indeed, the
elements $\tilde{w}$ of $G^0$ (and then of $G$) may be identified
(by solving the problem  (\ref{pbdirichlet})) with their traces $w$
on $\Gamma_{1}$. Moreover, as the functions $\tilde{w}$ are
harmonic, we may exhibit their local behavior in the vicinity of
any point $x_0 \in \Gamma _1$. Proceeding as in (\ref{approx}), (\ref{solmodifie2}) and taking only the decreasing exponential towards the domain (this is the classical approximation for the construction of a parametrix) we have:
\begin{equation}\label{fouriertangentiel}
\mathcal{F}(\tilde{w})(\xi _1, y_2)=\mathcal{F}(w)(\xi_1) e^{-|\xi
_1 |y_2},
\end{equation}
where $y_1, y_2$ are the tangent and the normal (inwards the
domain) vectors. Then, it is apparent that the $b$-energy is
concentrated in a layer close to $\Gamma _1$ and we may compute it
in an analogous way to the calculus that was done for the
$a$-energy (\ref{energiea3bis}). Indeed, using Parseval-Plancherel
Theorem and within our approximation, we have
\begin{eqnarray}\label{energieb}
b(\tilde{w}, \tilde{w})&=&\int _{-\infty }^{+ \infty } \D y_1 \int
_0 ^{+ \infty } \sum_{\alpha, \beta }^{}|\partial _{\alpha \beta }
\tilde{w} |^2 \D y_2\\
&=&  \int _{-\infty }^{+ \infty } \D \xi_1 \int _0 ^{+ \infty }
\Big(\xi _1 ^4 |\mathcal{F}(\tilde{w})|^2+ 2 \xi _1 ^2
|\mathcal{F}(\frac{ \partial \tilde{w}}{\partial y _2})|^2
+|\mathcal{F}(\frac{ \partial ^2\tilde{w}}{\partial y _2 ^2})|^2
\Big)\D y_2,\nonumber
\end{eqnarray}
hence, recalling (\ref{fouriertangentiel}) and integrating over
$y_2$, we get:
\begin{equation}\label{energieb2}
b(\tilde{w}, \tilde{w})=2 \int _{-\infty }^{+ \infty }| \xi _1 |
^3 |\mathcal{F}(w)|^2 \D \xi _1.
\end{equation}

Then, defining the pseudo-differential operator $Q(\frac{\partial
}{\partial s})$ of order $3/2$ with principal symbol
\begin{equation}\label{pseudodiffQ}
\sqrt{2} |\xi _1 |^{3/2},
\end{equation}
or equivalently as previously:
\begin{equation}\label{pseudodiffQ2}
\sqrt{2} (1+ |\xi _1 |^2)^{3/4},
\end{equation}
we have (always within our approximation):
\begin{equation}\label{energieb3}
b(\tilde{v}, \tilde{w})= \int _{\Gamma _1 }Q(\frac{\partial
}{\partial s})v  \ \overline{Q(\frac{\partial }{\partial s})w}\D
s.
\end{equation}

We observe that the operator $Q$ is only concerned with the trace
on $\Gamma_{1}$, so that we may either write $\tilde{v}$,
$\tilde{w}$ or $v$, $w$ in (\ref{energieb3}).

The formal asymptotic problem becomes:
\begin{equation}\label{pblaxasymptotic2}
\left\{\begin{array}{l}\textrm{Find }  \tilde{v}^{\varepsilon }\in G \textrm{ such that } \forall \tilde{w} \in G\\
\int _{\Gamma _0} P(\frac{\partial \tilde{v}^{\varepsilon }
}{\partial n}) \overline{P(\frac{\partial \tilde{w}  }{\partial
n}) }\D s + \varepsilon ^2 \int _{\Gamma _1} Q
(\tilde{v}^{\varepsilon })  \ \overline{Q(\tilde{w})}\D s =
\dualprod{ f}{w}.
\end{array}\right.
\end{equation}

\subsection{The formal asymptotics and its sensitive behaviour}
\noindent


In order to exhibit more clearly the unusual character of the
problem, we shall now write (\ref{pblaxasymptotic2}) under another
equivalent form involving only the traces on $\Gamma_{1}$.
Coming back to (\ref{pbdirichlet}), let us define $\mathcal{R}_0$
as follows. For a given $w \in C^{\infty }(\Gamma _1)$ we solve
(\ref{pbdirichlet}) and we take the trace of $\frac{\partial
\tilde{w}}{\partial n}$ on $\Gamma _0$, then
\begin{equation}\label{defdeR}
\frac{\partial \tilde{w}}{\partial n}_{|\Gamma _0} = \mathcal{R}_0
w.
\end{equation}

Using the regularity properties of the solution of
(\ref{pbdirichlet}), it follows that $\mathcal{R}_0 w$ is in
$C^{\infty }(\Gamma _0)$. In fact, $R_{0}$ is a smoothing
operator, sending any distribution into a $C^{\infty}$ function.
Then, (\ref{pblaxasymptotic2}) may be written as a problem for the
traces on $\Gamma _1$:
\begin{equation}\label{pblaxasymptotic3}
\left\{\begin{array}{l}\textrm{Find }  v^{\varepsilon }\in
H^{3/2}(\Gamma _1) \textrm{ such that } \forall w \in
H^{3/2}(\Gamma _1)\\
\int _{\Gamma _0} P(\frac{\partial   }{\partial s})
\mathcal{R}_0v^{\varepsilon } \overline{P(\frac{\partial
}{\partial s}) \mathcal{R}_0w }\D s + \varepsilon ^2 \int _{\Gamma
_1} Q(\frac{\partial   }{\partial s}) v^{\varepsilon }  \
\overline{Q(\frac{\partial }{\partial s})w}\D s = \int_{\Omega }
F\tilde{w} \D x,
\end{array}\right.
\end{equation}
where the configuration space is obviously $H^{3/2}(\Gamma _1)$.
The left hand side with $\varepsilon >0$ is continuous and
coercive.
 We then define the new operators
\begin{eqnarray}\label{operateur2}
\mathcal{A}&=& \mathcal{R}_0^*  P^* P \mathcal{R}_0\  \in
\mathcal{L}(H^s(\Gamma _1),H^{r }(\Gamma _0)), \forall s, r \in \R, \\
\mathcal{B}&=&   Q^* Q\  \in \mathcal{L}(H^{3/2}(\Gamma _1) , H^{-
3/2}(\Gamma _1)
\end{eqnarray}
where $\mathcal{R}_0^*  $ is the adjoint of $\mathcal{R}_0$ (which
is also smoothing)), (\ref{pblaxasymptotic3}) becomes
\begin{equation}\label{pblaxasymptotic4}
\Big( \mathcal{A} + \varepsilon ^{2}\mathcal{B}\Big)
v^{\varepsilon} =F, \textrm{ in } H^{- 3/2}(\Gamma _1).
\end{equation}
 Obviously, $\mathcal{B}$  is an elliptic
pseudo-differential operator of order 3, whereas $\mathcal{A}$ is
a smoothing (non local) operator.

 This problem is somewhat simpler than the initial one (as on a manifold of dimension 1), showing the interest of the formal asymptotics. It enters in a class of sensitive problems addressed in \cite{EgMeSa} section 2. It is apparent that the limit problem (for $\varepsilon =0$) has no solution in the distribution space for any $F$ not contained in $\mathcal{C}^\infty$. Indeed, on the compact manifold $\Gamma _0$, any distribution is in some $H^{-m} (\Gamma _0)$ space, which is send into $\mathcal{C}^\infty$ by the smoothing operator $\mathcal{A}$.

\begin{remark}\label{dernierrajout1}
The drastically non local character of the smoothing operator
$\mathcal{A}$ follows from the fact that it involves
$\mathcal{R}_0 $ and $ \mathcal{R}_0^*$ (see(\ref{defdeR})). This
is the reason why the problem may be reduced to another one on the
traces on $\Gamma_{1}$. The possibility of that reduction is a
consequence of our approximation, where the configuration space is
formed by harmonic functions.
\end{remark}

\end{document}